\documentclass[12pt]{article}

\usepackage[latin1]{inputenc}
\usepackage{verbatim}
\usepackage{exscale}
\usepackage{amssymb}
\usepackage{dsfont}
\usepackage{newlfont}
\usepackage{epsfig}
\usepackage{caption}
\usepackage{multirow}
\usepackage{pstricks}
\usepackage{pst-node}
\usepackage{pst-coil}
\usepackage[numbers]{natbib}  
\usepackage[breaklinks,colorlinks]{hyperref} 
\usepackage{amsmath} 
\usepackage{latexsym}

\newtheorem{Theorem}{\textit{Theorem}}
\newtheorem{Definition}{\textit{Definition}}
\newtheorem{Proposition}{\textit{Proposition}}
\newtheorem{Remark}{\textit{Remark}}

\title{Asymptotic for the quantization error for a Wiener process with Gaussian starting point}
\author{Luis A. Salomón\\{ \small Universidad de La Habana. Facultad de Matemática y Computación}\\ {\small email: \rm algoritmo.cu@gmail.com}}
\date{\footnotesize March, 2014}
\begin{document}

\maketitle

\begin{abstract}
The asymptotics for quantization error for a Wiener process with Gaussian starting point (GSP-Wiener process) is investigated. Using the classical methodology and some analytical approach a first result is obtained. We provide some further comments on the sharp asymptotic for the quantization error attending to a numerical experiment around the approximate solutions of the eigenvalues related with the covariance operator of the process.   
\end{abstract}



\vfill
\rule{7cm}{0.4pt}

\footnotesize
Accepted in Revista Ciencias Matemáticas
\normalsize
\hypersetup{citecolor=black}
\hypersetup{colorlinks, linkcolor=black}
\section{Introduction}

Functional quantization has been widely investigated from a practical and theoretical point of view.  The first approach to functional quantization is due to \citet{PagesLuschgy2002}. \ This theory is the natural extension to stochastic process of the optimal vector quantization of random vectors in finite dimension (usually $\mathbb{R}^d$), see for instance \citet{LuschgyGraf2000} for a deeper discussion.

This theory  studies the best approximation of a stochastic process in their path \ spaces by random vectors taking at most $n$ values. For the Gaussian process there is an important number of extensions around the rate of convergence of the quantization error, see for instance \citet{PagesLuschgy2004}, \citet{DereichLifshits2005}, \citet{PagesLuschgy2004a}, \citet{DereichFehringerMatoussiScheutzow2003} and \citet{DereichScheutzow2006}.

Let us consider a separable Hilbert space $(H,\langle\cdot,\cdot\rangle_H)$  with its natural $\sigma$-algebra. In this framework the most frequent choice for $H$ is $\mathcal{L}^2_T$, where $\mathcal{L}^2_T=\mathcal{L}^2([0,T],dt)$ with its usual norm. One considers a random variable $X$ (stochastic process) defined on a probability space $(\Omega, \mathcal{F},\mathbb{P})$ taking its values in $H$.

For $n\in\mathbb{N}$, the $n$-quantization problem for $X$ satisfying $\mathbb{E}||X||^2<+\infty$  consists in minimizing
\begin{equation}\label{eq: error}
\left(\mathbb{E}\min_{1\leq i\leq n}||X-a_i||_{\mathcal{L}^2_T}\right)^{1/2},
\end{equation}
over all the sets $\alpha=\{a_1,\cdots,a_N\}\subset H$ and $|\alpha|\leq n$, where the set $\alpha$  is called $n$-$codebook$ or $n$-quantizer. The minimal quantization error $e_n^2(X,\alpha)$ is defined  by
\begin{eqnarray*}
    e_n^2(X,\alpha)&=&\inf\left\{\mathbb{E}\min_{1\leq i\leq n}||X-a_i||_{\mathcal{L}^2_T}:\alpha\subset H,\;|\alpha|\leq n\right\}.
\end{eqnarray*}

We write $e_n^2(X)$ instead $e_n^2(X,\alpha)$ when not confusion can arise. From the $n$-quantizer we construct an approximation $\widehat{X}^{\alpha}:H\rightarrow \alpha\subset H$ of $X$, obtained by the rule of closest neighbor 
$$
\widehat{X}^{\alpha}=\pi_{\alpha}(X)=\sum_{i=1}^na_i\boldsymbol{1}_{C_{\alpha}(a_i)}(X),
$$
where $\{C_{\alpha}(a_i)\}_{1\leq i\leq n}$ is the Vorono\" {i} partition of $H$ defined as usual   
\begin{eqnarray*}
   C_{\alpha}(a_i)\subset V_{\alpha}(a_i) &\stackrel{\triangle}{=}& \{x\in
H:||x-a_i||_{\mathcal{L}^2_T}=\min_{1\leq j\leq n}||x-a_j||_{\mathcal{L}^2_T}\}.
\end{eqnarray*}

One of the most important property of optimal quantizers is the asymptotic behavior of the quantization error. It is easy to check that $e_n^2$ decreases to zero when the number of quantizers go to infinity. 

In the finite dimensional framework ($H=\mathbb{R}^d$) the asymptotic rate comes from the earlier work of \citet{Zador1982} and it was completely fulfilled in \citet{LuschgyGraf2000}.
 \begin{Theorem}\textbf{Rate of decay}\label{teo: zador}\\
Let $r>0$, assume that $\int_{\mathds{R}^d}|\xi|^{r + \eta}\mathbb{P}(d\xi)<+\infty$ for some $\eta>0$. Set $f\stackrel{\triangle}{=}d\mathbb{P}/d\lambda_d$, then

\begin{equation*}
 \lim_n\left(n^{\frac{r}{d}}\min_{(\mathds{R}^d)^{n}}e_{n,r}^r(X)\right)=J_{r,d}||f||_{\frac{d}{d+r}}<+\infty,
 \end{equation*}
 where
 \begin{eqnarray*}
    e_{n,r}(X)&=&\inf\left\{(\mathbb{E}\min_{a\in\alpha}||X-a||^r)^{1/r}:\alpha\subset \mathbb{R}^d,\;|\alpha|\leq n\right\},
 \end{eqnarray*}
 and
  $$ 
  ||g||_p\stackrel{\triangle}{=}\left(\int_{\mathds{R}^d}|g(\xi)|^{p}\mathbb{P}(d\xi)\right)^{1/p}
  \quad\text{for every}\; p\in(0,+\infty).
  $$
 \end{Theorem}
The positive real constant $J_{r,d}$ corresponds to the uniform distribution on $[0,1]^d$. One knows that $J_{r,1}=1/(2^r(r+1))$, $J_{2,2}=5/(18\sqrt{3})$. When  $d\geq 3$, $J_{r,d}$ is unknown. However, when  $d\rightarrow +\infty$ the following asymptotic expansion holds $J_{r,d}=d/(2\pi e)^{r/2}+o(d)$ .

In the functional quantization settings  a similar results exist for a wide class of Gaussian stochastic process. For the Wiener process $W$ in $[0,T]$ in the quadratic case (see \citet{PagesLuschgy2004} for details) it is known that
 \begin{equation*}\label{eq: e_n Wiener}
 e_{n}(W)\sim \frac{\sqrt{2}T}{\pi} (\ln n)^{-1/2}.
 \end{equation*}

This paper is intended as an attempt to derive some asymptotics for the quantization error of a GSP-Wiener process. This Gaussian process could be described  in $[0,T]$ as $Z_k+W_t$, where $Z_k\sim N(0,k)$ and $W$ is a Wiener process, $Z_k$ is independent with $W$. Under this definition it seems that a GSP-Wiener  should have the same asymptotic rate of the Wiener process. The paper is organized as follows: next section contains a brief summary on the fundamental results around the convergence rate of the quantization error, section 3 is devoted to the main result of the paper and finally some comments and further remarks are provided at the end of this work.
 
\section{Convergence rate in functional quantization}

Some properties of the finite dimensional case hold in the Hilbert space framework $H=\mathcal{L}^2_T$. Several authors have studied the convergence rate of the functional quantization of $e_n^2$. The best general reference here is due to \citet{PagesLuschgy2002}, \citet{PagesLuschgy2004}, \citet{Dereich2003} and references therein. In the infinite dimensional Gaussian setting exist a general version of Theorem \ref{teo: zador}.

\begin{Theorem}\textbf{Rate of Decay}\label{teo: zador Hilbert}\\
If $X$ is a centered $H$-valued random vector with a Gaussian distribution $\mathbb{P}$, covariance operator $\Gamma_C$  and $(u^X_j)_{j\geq1}$ is any orthonormal subset of $H$, such as the reproducing kernel Hilbert space $K_X$ satisfies $K_X\subset\mathrm{cl}\,\mathrm{span}\{u^X_j,j\in\mathbb{N}\}$. Let us define
\begin{equation*}
\mu_j=\mathrm{Var}\langle u^X_j,X\rangle=\langle u^X_j,\Gamma_C(u^X_j)\rangle\quad\text{and}\quad\Sigma_m=(\langle u^X_j,\Gamma_C(u^X_k)\rangle)_{0\leq j,k \leq m},
\end{equation*}
and for $n\in \mathds{N}$ set
\[
g_n(m)=e_n(\mathcal{N}(0,\Sigma_m)).
\]
If we assume that $\{u^X_j,j\in\mathbb{N}\}\subset\mathrm{cl} (K_X)$, then $\mathrm{det}\,\Sigma_m>0$ and
\[
\lim_{n\rightarrow \infty}n^{1/m}g_n(m)=Q(m)\quad \text{for every} \quad m\geq 1,
\]
where $Q(m)\in (0,\,+\infty)$ and
\[
Q(m)\thicksim\left(m(\mathrm{det}\,\Sigma_m)^{1/m}\right)^{1/2}\quad\text{as}\quad m\rightarrow \infty.
\]
In particular
\[
\lim_{m\rightarrow\infty}Q(m)=0.
\]
\end{Theorem}
The previous theorem is quite general. More specific asymptotics for the quantization error exist for a large variety of Gaussian process. The methods used to establish the rate of convergence rely on the behavior of the eigenvalues of the covariance operator, varying functions, small ball probabilities and Shannon-Kolmogorov $\epsilon$-entropy.

The first approach to the asymptotics of quantization error is due to \citet{PagesLuschgy2002}. They give us upper and lower bounds for $e_n$ using  the eigenvalues behavior (varying function) of the covariance operator  and Shannon-Kolmogorov's $\epsilon$-entropy respectively. (See Theorem 4.6 and 4.12 in \citet{PagesLuschgy2002}).

In \citet{PagesGrafLuschgy2003} and \citet{Dereich2003} a different approach is proposed using varying functions, small ball probabilities and their link with the quantization error, see for instance Theorem 1.2 in \citet{PagesGrafLuschgy2003} and Theorem 2.3 in \citet{Dereich2003}. Sharp bounds are given for a wide class of Gaussian processes. (See also \citet{DereichFehringerMatoussiScheutzow2003}).

The most important results for sharp asymptotics of quantization error are detailed in \citet{PagesLuschgy2004,PagesLuschgy2004a}. The general idea is mainly based in regular varying functions and Shannon-Kolmogorov's $\epsilon$-entropy.  We recall here the Theorem 2.2 (a) in  \citet{PagesLuschgy2004} which provides the sharp rate of convergence for the distortion.
\begin{Theorem}\textbf{ }\label{teo: 2.2 Pages 2004}\\
If $X$ is a Gaussian process with eigenvalues $\lambda_j\sim\varphi(j)$ as $j\rightarrow\infty$, where $\varphi:(s,\infty)\rightarrow(0,\infty)$ is a decreasing, regularly varying function at infinity of index $-b<-1$ for some $s\geq0$. Set, for every $x>s$,
$$
\psi(x)\stackrel{\triangle}{=}\frac{1}{x\varphi(x)},
$$
then
\begin{equation}\label{eq: equi e_n Pages2004}
    e_n(X)\sim \left(\left(\frac{b}{2}\right)^{b-1}\frac{b}{b-1}\right)^{1/2}\psi(\ln n)^{-1/2}.
\end{equation}
\end{Theorem}
The most prevalent form for $\varphi$ is
$$
\varphi(x)=cx^{-b}(\ln x)^{-a},\quad b>1,\;a\in\mathds{R},\;x>\max\{1,e^{-a/b}\},
$$
and (\ref{eq: equi e_n Pages2004}) turns to
$$
e_n(X)\sim \left(c\left(\frac{b}{2}\right)^{b-1}\frac{b}{b-1}\right)^{1/2}(\ln n)^{-(b-1)/2}(\ln\ln n)^{-a/2}.
$$

For the particular case of diffusions the Shannon-Kolmogorov's $\epsilon$-entropy plays a fundamental role, see for instance \citet{PagesLuschgy2006}, \citet{DereichScheutzow2006} and \citet{Dereich2008a}. Some asymptotics are also provided for the $d$-dimensional Brownian motion.  In \citet{Salomon2011} a brief summary of asymptotics for the quantization error is presented for a wide class of stochastic processes.

Theorem \ref*{teo: 2.2 Pages 2004} establishes the eigenvalues importance for the rate of the quantization error. The theorems which states the asymptotic of eigenvalues comes from the works of Widom and Rosenblatt (see Theorem 1 in \citet{Widom1963} and Theorem 3 in \citet{Rosenblatt1963} for a deeper discussion).
\section{Main Results}
Let us denote by $Z^{(k)}$ a GSP-Wiener process. It is obvious that $Z^{(k)}$ is a Gaussian process and has similar properties to those of  Brownian motion. The formal definition could be written as follows:
\begin{Definition}\textbf{ }\\ \label{def: Wiener truncated}
A Gaussian stochastic process  $\{Z_t^{(k)},\,t\in [0,T]\}$ for $k\geq0$ is a GSP-Wiener process if
$$
Z_t^{(k)}=W_{k+t},\quad\forall t\in [0,T],
$$
where $W$ is a Wiener process.
\end{Definition}
This process satisfies the following conditions
\begin{Proposition}\textbf{ }\\ \label{prop: Wiener}
Let $Z^{(k)}$ be a GSP-Wiener process hence 
\begin{enumerate}
\item $Z_0^{(k)}\sim N(0,k)$.
\item The process $Z^{(k)}$ has stationary and independent increments.
\item $\mathbb{E}(Z_t^{(k)})=0$ and $\mathrm{cov}(Z_t^{(k)},Z_s^{(k)})=(k+t)\wedge (k+s)$.
\item $Z^{(k)}$ is not self-similar.
\end{enumerate}
\end{Proposition}
The proof of previous proposition is a simple consequence of  $Z^{(k)}$ definition and classical properties of the Wiener process.

\subsection{Asymptotic for the quantization error}
The main result for the stochastic process $Z^{(k)}$ read as follows
\begin{Theorem}\textbf{ }\\ \label{teo: e_n truncated Wiener}
If $Z^{(k)}$ is a \ stochastic process as in \ Definition \ref{def: Wiener truncated} then the eigenvalues associated to the process $(\lambda_\ell)_{\ell\geq1}$ satisfies that $\lambda_\ell\approx \ell^{-2}$ and
\begin{equation}\label{eq: e_n Z approx}
e_n(Z^{(k)})=\Theta\left((\ln n)^{-\frac{1}{2}}\right) ,
\end{equation}
where $\Theta(f(n))=O(f(n))\cap \Omega(f(n))$.
\end{Theorem}
\subsubsection*{Proof}
The underlying idea here is to use Theorem \ref{eq: equi e_n Pages2004}. Hence the first step is to obtain the eigenvalues associated to the process.

The covariance function for $Z^{(k)}$ is defined by
$$
K_{Z^{(k)}}(t,s)=\mathrm{cov}(Z_t^{(k)},Z_s^{(k)})=\mathbb{E}(Z_t^{(k)}\cdot Z_s^{(k)})=(k+t)\wedge (k+s).
$$
The classical approach to obtain the eigenvalues comes from solution of the integral equation:
$$
\int_{0}^{T}K_{Z^{(k)}}(t,s)\varphi(s)ds=\lambda \varphi(t),\qquad\forall t\in [0,T].
$$
Working in the previous equation it follows that
\begin{eqnarray}\label{eq: varphi despeje}
    \lambda \varphi(t)&=&\int_{0}^{T}K_{Z^{(k)}}(t,s)\varphi(s)ds=\int_{0}^{T}((k+t)\wedge (k+s))\varphi(s)ds \nonumber\\
    &=&\int_{0}^{t}s\varphi(s)ds+t\int_{t}^{T}\varphi(s)ds+k\int_{0}^{T}\varphi(s)ds.
\end{eqnarray}
Differentiating twice with respect to $t$ yields to the following system:
\begin{eqnarray*}
    \lambda \varphi'(t)&=&\int_{t}^{T}s\varphi(s)ds\\
    \lambda \varphi''(t)&=& -\varphi(t).
\end{eqnarray*}
It is straightforward that a general solution of previous equation has the form
\begin{equation}\label{eq: sol EDO KL}
\varphi(t)=\alpha \sin\left(\frac{T}{\sqrt{\lambda}}\right)+\beta \cos\left(\frac{T}{\sqrt{\lambda}}\right).
\end{equation}
At this step we need to find the values for $\alpha$ and $\beta$. Using equation (\ref{eq: varphi despeje}) and equation (\ref{eq: sol EDO KL}) for $t=0$ we obtain that 
\begin{equation*}\label{eq: beta}
\beta=\frac{k}{\lambda}\int_{0}^{T}\varphi(s)ds.
\end{equation*}
Differentiating equation (\ref{eq: sol EDO KL}) for $t=0$ it follows that
\begin{equation*}\label{eq: beta}
\alpha=\lambda^{-\frac{1}{2}}\int_{0}^{T}\varphi(s)ds.
\end{equation*}
It is obvious that we have determined the values for $\alpha$ and $\beta$; however, both cases depend on the unknown $\varphi$. Hence it seems that we are in a similar problem. However using the fact that $\lambda \varphi'(T)=0$ and the expression for $\varphi'$ obtained by equation (\ref{eq: sol EDO KL}) it follows that
\begin{eqnarray*}
    0&=&\varphi'(T)=\frac{\alpha}{\sqrt{\lambda}} \cos\left(\frac{T}{\sqrt{\lambda}}\right)- \frac{\beta}{\sqrt{\lambda}}\sin\left(\frac{T}{\sqrt{\lambda}}\right)\\
    &=& \lambda^{-1}\int_{0}^{T}\varphi(s)ds\left( \cos\left(\frac{T}{\sqrt{\lambda}}\right)- \frac{k}{\sqrt{\lambda}}\sin\left(\frac{T}{\sqrt{\lambda}}\right)\right)\\
    &=&\cos\left(\frac{T}{\sqrt{\lambda}}\right)- \frac{k}{\sqrt{\lambda}}\sin\left(\frac{T}{\sqrt{\lambda}}\right).
\end{eqnarray*} 
Hence we have
\begin{equation}\label{eq: sol cot}
\cot\left(\frac{T}{\sqrt{\lambda}}\right)=\frac{k}{\sqrt{\lambda}}.
\end{equation}
If we write $x=\frac{T}{\sqrt{\lambda}}$ then previous equation reads as
\begin{equation}\label{eq: sol cot2}
\cot\left(x\right)=\frac{k}{T}x.
\end{equation}
Hence, solving this last equation should give us the corresponding solution for $\lambda$. Let us analyze this equation. The $\cot(x)$ function is decreasing in each interval $((\ell-1)\pi,\ell\pi)$, for all $\ell\in\mathbb{N}$ and $\frac{k}{T}x$ is an increasing function. Therefore in each $((\ell-1)\pi,\ell\pi)$ there exists one and only one $x^{(k)}_\ell$ for which $\cot\left(x^{(k)}_\ell\right)=\frac{k}{T}x^{(k)}_\ell$. Knowing that $k\geq0,T>0$ it is easy to check that the mentioned $x^{(k)}_\ell$ is actually in $((\ell-1)\pi,(2\ell-1)\frac{\pi}{2}]$. Therefore the solution for equation (\ref{eq: sol cot}) is given by the sequence
$$
(\lambda^{(k)}_\ell)_{\ell\in\mathbb{N}}=\left(\left(\frac{T}{x^{(k)}_\ell}\right)^2\right)_{\ell\in\mathbb{N}}.
$$
By the same argument used before is straightforward that
$$
\lambda^{(k)}_\ell\in \Bigg[\frac{T^2}{((2\ell-1)\frac{\pi}{2})^2},\frac{T^2}{((\ell-1)\pi)^2} \Bigg),\quad\forall\ell\in\mathbb{N}.
$$
Let us write then
$$
\lambda^{(k)}_\ell=c^{(k)}_\ell\frac{T^2}{((2\ell-1)\frac{\pi}{2})^2},
$$
where $(c^{(k)}_\ell)_{\ell\in\mathbb{N}}$ is a sequence that should satisfies $1\leq c^{(k)}_\ell<3/2$, for all $k$ and $\ell\geq2$, and also for $k=0$, $c^{(0)}_\ell=1$, $\forall\ell\in\mathbb{N}$. For the case $k=0$ is quite obvious that $Z^{(0)}$ is a Wiener process $[0,T]$ and we retrieve its classical eigenvalues. In fact we can easy check that $\lambda^{(k)}_\ell\approx \lambda^{W}_\ell$, for all $\ell$, hence by Lemma 4.11 in \citet{PagesLuschgy2002} it follows that
$$
e_n(Z^{(k)})\approx e_n(W).
$$
In fact for $n$ large enough
$$
\sqrt{2}\frac{T}{\pi}(\ln n)^{-1/2} \leq e_n(Z^{(k)}) \leq \sqrt{3}\frac{T}{\pi}(\ln n)^{-1/2},
$$
and equation (\ref{eq: e_n Z approx}) follows after that and the proof is complete
\begin{flushright}
$\blacksquare$
\end{flushright}
\begin{Remark}\textbf{ }\\ \label{rem: e_n truncated Wiener}
If the sequence $(c^{(k)}_\ell)_{\ell\in\mathbb{N}}$ is convergent, $i.e.$,
$$ c^{(k)}_\ell\rightarrow c^{(k)}_\infty,\quad \text{when} \;\ell\rightarrow \infty,$$
then  we have a sharp asymptotics for $e_n(Z^{(k)})$
$$
e_n(Z^{(k)})\sim \sqrt{2c^{(k)}_\infty}\frac{T}{\pi}(\ln n)^{-1/2}.
$$
\end{Remark}
\subsection{Some numerical calculations}
As was shown in previous subsection the eigenvalues $(\lambda^{(k)}_\ell)_{\ell\in\mathbb{N}}$ for the process $Z^{(k)}$ are the solution of equation 
$$
\cot\left(\frac{T}{\sqrt{\lambda^{(k)}_\ell}}\right)=\frac{k}{\sqrt{\lambda^{(k)}_\ell}},\qquad \lambda^{(k)}_\ell\in \Bigg[\frac{T^2}{((2\ell-1)\frac{\pi}{2})^2},\frac{T^2}{((\ell-1)\pi)^2} \Bigg),\quad\forall\ell\in\mathbb{N}
$$
The solution for previous equation could be only solved by numerically. There are several known root finding method as bisection, Newton among others to solve this problem. We use  classical Newton's method to obtain these solutions.

In order to study the behavior of these solutions we compute the first 1000 eigenvalues for the process $Z^{(k)}$ in $[0,1]$ for 100 different values of $k$ ($k_i=(i-1)\cdot 10^{-2}$, where $i=1,2,\cdots,100$). We work with equation (\ref{eq: sol cot2}).

The implementation of  Newton's method requires a starting point $x^{(0)}_\ell$. We take this $x^{(0)}_\ell=10^{-5}+ (\ell-1)\cdot\pi$ on each interval $((\ell-1)\pi,(2\ell-1)\frac{\pi}{2}]$ for $\ell=1,2,\cdots,1000$. The choice for $x^{(0)}_\ell$ was determined  numerically for some fixed $k^\star\in [0,1)$. This first study shows in some intervals that for greater values of $x^{(0)}_\ell$  the solution given by Newton's method was the same for two intervals. However for values of $x^{(0)}_\ell$ near to $(\ell-1)\pi$ we find the correct solution in each one, at least for the first 1000 eigenvalues.

\begin{figure}[!t]
\begin{center}
\includegraphics[width=12 cm, height=6 cm]{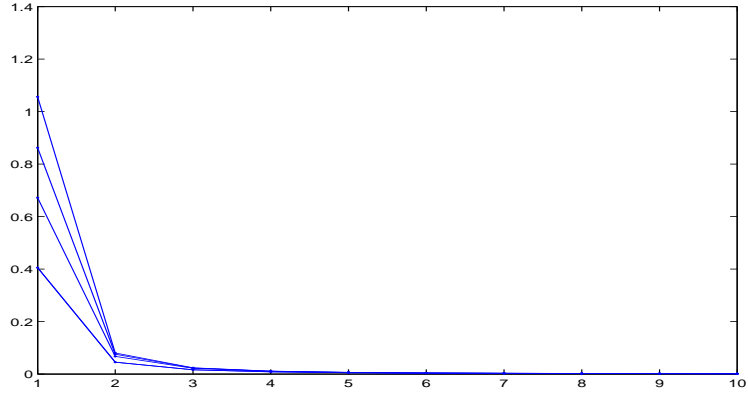}
\end{center}
\caption{Behavior of the first 10 eigenvalues for $Z^{(0)}$, $Z^{(0.3)}$, $Z^{(0.5)}$ and $Z^{(0.7)}$ from bottom to top.}\label{fig: EigenZ4}
\end{figure}

In Figure \ref{fig: EigenZ4} it is shown the behavior of 10 eigenvalues for four  stochastic processes $Z^{(k)}$. It is easy to see that in all cases, it exhibits a similar behavior. The most important difference relies on the initial $\lambda^{(k)}_1$. 

\begin{figure}[!htb]
\begin{center}
\includegraphics[width=12 cm, height=6 cm]{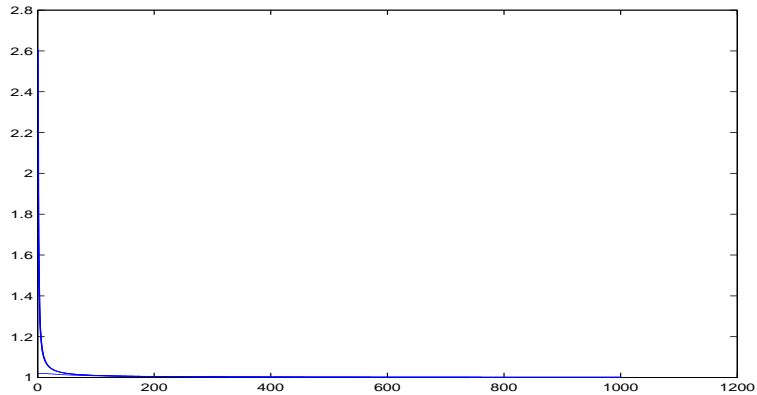}
\end{center}
\caption{Behavior of $\widehat{c}^{(k)}_\ell$ for  $Z^{(0.3)}$, $Z^{(0.5)}$ and $Z^{(0.7)}$.}\label{fig: EigenC}
\end{figure}

In Figure \ref{fig: EigenC} we can observe for the processes $Z^{(0.3)}$, $Z^{(0.5)}$ and $Z^{(0.7)}$ the behavior of
$$
\widehat{c}^{(k)}_\ell=\widehat{\lambda}^{(k)}_\ell/\widehat{\lambda}^{(0)}_\ell,\; \forall\ell=1,2,\cdots,1000. 
$$

The graph for the three processes are almost the same, however it is easy to see that $\widehat{c}^{(k)}_\ell$ shows, in all cases, a convergent behavior to 1. Thus Remark \ref{rem: e_n truncated Wiener} appears to be true. It seems, by this numerical experiment, that
$$ c^{(k)}_\ell\rightarrow 1,\quad \text{when} \;\ell\rightarrow \infty.$$
\section{Conclusions}
This work presents an asymptotic result for an specific Gaussian process  $Z^{(k)}$. As expected the quantization error for $Z^{(k)}$ exhibits a similar behavior as the Wiener process itself. This result could be used with theoretical purposes. The most important point to note here is the numerical results obtained for the eigenvalues. It seems that a further analysis on these numerical computations of the eigenvalues for  $Z^{(k)}$ could serve to find the specific behavior of the sequence $(c^{(k)}_\ell)_{\ell\in\mathbb{N}}$.

\bibliographystyle{plainnat}
\addcontentsline{toc}{section}{\textbf{References}}
\bibliography{BIBLIOGRAFIA_DOCTORAT}
\end{document}